\documentclass[a4paper,11pt]{amsart}
\usepackage{amsmath, amssymb, amsfonts, amsthm, enumerate, appendix, pdfpages}

\newtheorem{mythm}{Theorem}
\newtheorem{mylemma}{Lemma}
\newtheorem{myprop}{Proposition}

\def\Q{\mathbb Q}
\def\C{\mathbb C}
\def\Z{\mathbb Z}
\def\O3{$O_3$}
\def\SO3{$SO_3$}
\def\A4{$ A_4 $}
\def\S4{ $ S_4 $}
\def\Al5{ $ A_5 $}
\def\GL{\mathrm{GL}}

\begin{document}

\title[Self-dual Artin representations of dimension three]{Self-dual Artin representations of dimension three \\ \tiny{ (with an appendix by David E. Rohrlich)}}
\author{Aditya Karnataki}
\date{}

\begin{abstract}
We give an unconditional proof that self-dual Artin representations of $\Q$ of dimension $3$ have density 0 among all Artin representations of $\Q$ of dimension $3$. Previously this was known under the assumption of Malle's Conjecture.

\end{abstract}

\maketitle

\tableofcontents

\begin{section}{Introduction}
It is expected that essentially self-dual motives (i.e. motives that are dual to a Tate twist of themselves) should occur with density $0$. This statement has not yet been formulated for motives with weights higher than $0$, but for Artin motives of fixed dimension, this is a precise question as shown in \cite{ABCZ}. Let $F$ be a number field, and if $\rho$ is an Artin representation of $F$, let $q(\rho)$ be the absolute norm of the conductor ideal $\mathfrak{q}(\rho)$. We denote by  $ \vartheta_{F, n}(x) $ the number of isomorphism classes of Artin representations over $F$ of dimension $n$ with $q(\rho) \le x$ and by $ \vartheta^{sd}_{F, n}(x) $ the number of isomorphism classes of self-dual Artin representations over $F$ of dimension $n$ with $q(\rho) \le x$. Rohrlich proved in \cite{Rohrlich} that the quotient $\vartheta^{sd}_{\Q, 2}(x) / \vartheta_{\Q, 2}(x)$ goes to $0$ as $x$ goes to $\infty$. Thus, our density $0$ expectations are true for dimension $2$. (The $1$-dimensional case is elementary.) He proved the same result for $\Q$ and dimension $3$ under a weak form of Malle's Conjecture. In this paper, we remove this condition, viz. we confirm unconditionally,

\begin{mythm} 
\label{MainThm}

$$ \lim_{x \to \infty} \frac{ \vartheta^{sd}_{\Q, 3}(x)} { \vartheta_{\Q, 3}(x)}  = 0. $$

\end{mythm} We can replace $\Q$ by any number field and ask a similar question. But that case seems considerably harder. In dimension $1$, the density result for a general number field follows from work of Taylor \cite{T}.

Before describing our work, we set some notations. For a finite extension $K/F$ of number fields, we denote by $\mathfrak{d}_{K/F}$ the relative discriminant ideal and by $d_{K/F}$ its absolute norm. For $F = \Q$, we simply write $\mathfrak{d}_K$ and $d_K$. We denote by $\eta_{F, m}(x)$ the number of extensions $K/F$ inside a fixed algebraic closure $\bar{F}$ such that $[K:F] = m$ and $d_{K/F} \le x$. Also, if $T$ is a transitive subgroup of the symmetric group $S_m$, we denote by $\eta^{T}_{F, m}(x)$ the number of extensions $K/F$ for which $Gal(L/F) \cong T$ as permutation groups, where $L/F$ is the normal closure of $K$ over $F$ and $Gal(L/F)$ is viewed as a permutation group via its action on conjugates of a primitive element of $K$ over $F$.

We now describe the structure of this paper. In section $2$, we shall use results from \cite{Rohrlich} to see that the problem reduces to bounding the number of {\it irreducible} self-dual Artin representations of $\Q$. 

Thus, we wish to count irreducible self-dual Artin representations. Such a representation has to be either orthogonal or symplectic. Since the dimension is odd, we are reduced to the orthogonal case. We are thus reduced to analyzing irreducible orthogonal finite subgroups $G$ of $\GL_3(\C)$, where $G = \rho(Gal(K / \Q))$ and $K$ is the fixed field of ker $\rho$. 

Our general strategy is to replace $\vartheta$, which counts conductors of representations, by $\eta$, which counts discriminants of number fields. This method has been used in earlier works such as Bhargava-Ghate (\cite{BG}), Kl\"{u}ners (\cite{K}). We then appeal to results of Bhargava (\cite{Bhargava01}, \cite{Bhargava02}) and Bhargava, Cojocaru, and Thorne (\cite{BCT}).

We shall divide our analysis into two cases : ($1$) Those $G$ which are contained in \SO3 and, ($2$) those which are not. In section $3$, we analyze the subgroups occurring in case ($1$). Using bounds on ramification of primes, we obtain bounds in this case in terms of $\eta^{A_4}_{\Q, 4}(x)$, $\eta^{S_4}_{\Q, 4}(x)$, and $\eta^{A_5}_{\Q, 5}(x)$. 

Having obtained these, we turn our attention to case ($2$). We further divide this case into two parts. Part $1$ is the case where $-1 \not \in G$. We analyze this case in section $5$. If $-1 \not \in G$, then we show that $G \cong S_4$ and $\rho$ is monomial. Here, monomial is equivalent to being induced from a character of a subgroup of finite index. We remark that this $\rho$ is the so-called standard representation of $S_4$, which arises as a 3-dimensional irreducible summand of the permutation action of $S_4$ on an ordered basis of a 4-dimensional complex vector space. (See \cite{FH}, page $18$, for more details.) It can be seen that this representation is induced from a quadratic character of a subgroup of index three, which corresponds to a quadratic character of a cubic subextension. Therefore, we reduce our problem to counting such cubic extensions and quadratic characters, using the interplay between conductors and discriminants, and obtain bounds in terms of $\eta_{\Q, 3}(x)$. Then we can apply asymptotics for $\eta_{\Q, 3}(x)$, which is the main result of \cite{DH}, to get required bounds.

In section $6$, we deal with the case where $-1 \in G$. This implies that $G$ can be written as $H \times \{ \pm 1 \}$. In this case, in addition to monomial representations coming from $H \cong A_4$ or $S_4$, we must also contend with primitive representations coming from $H \cong A_5$. Here, being primitive is equivalent to not being induced from any proper subgroup. It is worth noting that the irreducible primitive case corresponds to the dominant term in all our analysis and it is only the power-saving result of Bhargava, Cojocaru and Thorne that helps us establish our result.

Finally, in section $7$, we combine results from sections $4$, $5$, $6$ to get the main theorem. 

\end{section}

\begin{section}{Acknowledgments}

I am deeply grateful to David Rohrlich for directing my attention towards this problem and for providing me with helpful hints. I would also like to thank Shaunak Deo for his indispensable inputs and all the innumerable discussions. I would like to thank Dipendra Prasad for his guidance. I am also thankful to the anonymous referee for making valuable comments and suggestions.

\end{section}

\begin{section}{Reduction To The Irreducible Case}

As in \cite{Rohrlich}, we have \begin{equation} \label{eq:selfdual}   \vartheta^{sd}_{\Q, 3}(x) =  \vartheta^{ab, sd}_{\Q, 3}(x) +  \vartheta^{1+2, sd}_{\Q, 3}(x) +  \vartheta^{irr, sd}_{\Q, 3}(x), \end{equation} where  

(a) $ \vartheta^{ab, sd}_{\Q, 3}(x) $ is the number of abelian self-dual Artin representations of $\Q$ of dimension $3$ with $q(\rho) \le x$,

(b) $ \vartheta^{1+ 2, sd}_{\Q, 3}(x) $ is the number of isomorphism classes of self-dual Artin representations of $\Q$ of dimension $3$ of the form $ \rho \cong \rho' \oplus \rho'' $ with $\rho'$ one-dimensional, $\rho''$ irreducible and two-dimensional and $q(\rho')q(\rho'') \le x$, and

(c) $ \vartheta^{irr, sd}_{\Q, 3}(x) $ is the number of irreducible self-dual Artin representations of $\Q$ of dimension $3$ with $q(\rho) \le x$.

From Theorem 2 of \cite{Rohrlich}, we see that $$ \vartheta^{ab, sd}_{\Q, 3}(x) = O(x( \mbox{log }x)^2), $$ while from equation (80) of the same paper, we see that $$ \vartheta^{1+ 2, sd}_{\Q, 3}(x) \ll x^{2 - \epsilon}. $$ Since, by Theorem 1 of \cite{Rohrlich}, $$  \vartheta^{ab}_{\Q, 3}(x) \sim O(x^2(\mbox{log }x)^2), $$ we see that, if we can prove \begin{equation} \label{eq:main} \vartheta^{irr, sd}_{\Q, 3}(x) = O(x^{2 - \epsilon}), \end{equation}then we conclude that the self-dual representations have density zero.
\end{section}

\begin{section}{Finite Irreducible Orthogonal Subgroups of $\GL_{3}(\C)$}
We are interested in irreducible self-dual Artin representations. By definition of self-dual representations, these have to be either orthogonal or symplectic, i.e. their image is contained either in $O_n(\mathbb{R})$ or $Sp_{2n}(\mathbb{C})$. (In particular the trace has to be real for these representations.) But since the dimension is odd, these have to be in \O3, where \O3 is the orthogonal group of real $3 \times 3$ matrices. We shall first concentrate on finite subgroups of \SO3. Referring to \cite{Artin}, chapter $5$, we see that every finite subgroup $G$ of \SO3 is one of the following :

\begin{enumerate}
\item $C_k$ : The cyclic group of rotations by multiples of $2 \pi / k$ about a line
\item $D_k$ : The dihedral group of symmetries of a regular $k$-gon
\item \A4 : The alternating group on $4$ variables 
\item \S4 : The symmetric group on $4$ variables 
\item \Al5 : The alternating group on $5$ variables 

\end{enumerate}

The cyclic groups and dihedral groups do not possess irreducible representations of dimension three. The last three groups do have irreducible $3$-dimensional representations.  Note that\S4 has two irreducible $3$-dimensional representations, but the image of only one of them is contained in \SO3. We call a subgroup $G$ of $\GL_n(\C)$ \emph{irreducible} if the inclusion $ i : G \rightarrow \GL_n(\C) $ is an irreducible representation of $G$.

\end{section}

\begin{section}{A Bound on Discriminants for finite subgroups of \SO3}
In this section, $\rho$ is an irreducible self-dual Artin representation and $K$ denotes the fixed field of ker $\rho$. We know that $\rho(Gal(K / \Q))$ is a finite irreducible subgroup of $O_3$. We divide our analysis into two cases, depending upon whether $\rho(Gal(K / \Q))$ is a subgroup of $SO_3$ or not. Thus, we write : \begin{equation} \label{eq:division} \vartheta^{irr, sd}_{\Q, 3}(x) = \vartheta_1(x) + \vartheta_2(x), \end{equation} where $$ \vartheta_1(x) = \sum_{ \substack{\rho \\ \rho(Gal(K / \Q)) \subset  SO(3) \\ q(\rho) \le x }} 1 $$ and $$ \vartheta_2(x) = \sum_{ \substack{\rho \\ \rho(Gal(K / \Q)) \not \subset  SO(3) \\ q(\rho) \le x }} 1. $$
For the rest of the section, we focus on bounding $\vartheta_1(x)$. From here on, we identify $Gal(K / \Q)$ with its image under $\rho$. Thus, we assume $Gal(K / \Q) \subset SO_3$. We have seen in the previous section that this implies that $Gal(K / \Q) $ is isomorphic to $ A_4, S_4, \text{ or } A_5$, and we write $m$ for the degree of the permutation group in question. Thus, $m = 4$ in the first two cases and $m=5$ for the third. In all that follows, $M$ is any subfield of $K$ with $[M: \Q] = m$. The choice of $M$ is arbitrary, but the normal closure of $M$ is $K$ for every one of them.

\begin{myprop} \label{SO3A4}
If $Gal(K/\Q) \cong A_4$, $$ d_M \le c q(\rho)^{3/2} $$ with an absolute constant $ c > 1$.
\end{myprop}

\begin{proof}
We quote a standard bound ({\it cf.} \cite{Serre}, p. 127, Proposition 2), which is \begin{equation} \label{eq:ram} d_M \le c \prod_{\substack{p | d_M \\ p > m}} p^3 \end{equation} with $ c = 2^{11}3^{7}$.

Now, if $ p > 4 $ and $p | d_M$, then $\rho$ restricted to the inertia group $I$ for any prime $\mathfrak{p}$ above $p$ factors through its tame quotient (since $2$ or $3$ are the only wildly ramified primes for im $\rho \cong $ \A4) and hence, by the formula for (local) Artin conductor, \begin{equation}\label{eq:conductor} ord_p(q(\rho)) = \dim(V/V^I) \end{equation} where $V$ is the space of $\rho$ and $V^I$ is the subspace of inertial invariants.

{\bf Case 1: I is a cyclic subgroup of order 2.}
Since all elements of order $2$ are conjugate to each other in \A4, only one computation will suffice for all the three subgroups of order $2$. We see from a character table (see e.g. \cite{FH}) and Frobenius reciprocity that

$$ \mbox{Multiplicity of trivial character in } \rho|_{I} = \frac{3(1) + (-1)(1) }{2} = 1. $$ Hence, $ \dim(V^I) = 1 $. So that $\dim(V/V^I) = 2$.

Alternatively, we can also argue without referring to a character table as follows : We see that since the determinant of $\rho$ is $1$, the image under $\rho$ of a non-trivial element is conjugate to 

$$  \left( \begin{array}{ccc}
1 &   &   \\
  & -1 &  \\
  &   & -1 \end{array} \right)  $$ from which it is immediate that $ \dim(V^I) = 1 $. We record this method here, since it is this method that will be useful in the latter sections.

{\bf Case 2: I is a cyclic subgroup of order 3.} There are two conjugacy classes, each containing $4$ elements of order $3$, which cover all the elements of order $3$ in \A4. The character of our $3$-dimensional representation is valued $0$ on both of these classes. Hence, for restriction to any subgroup of order $3$, we see that $$ \mbox{Multiplicity of trivial character} = \frac{3(1) + 0(1) + 0(1)}{3} = 1. $$ Hence, $ \dim(V^I) = 1 $. So that $ \dim(V/V^I) = 2$.

Therefore, by (\ref{eq:conductor}), we see that $$ q(\rho) \ge \prod_{\substack{p | q(\rho) \\ p > m}} p^2 $$ and combined with (\ref{eq:ram}), it completes the proof.

Again, we can argue without using a character table : Since the determinant of $\rho$ is $1$, the image under $\rho$ of a non-trivial element is conjugate to

$$ 
\left( \begin{array}{ccc}
1 &   &   \\
  & \zeta &  \\
  &   & \zeta^2 \end{array} \right) $$ where $\zeta$ is a primary cube root of unity. Thus, it is immediate that $ \dim(V^I) = 1 $.

\end{proof}

\begin{myprop} \label{SO3S4}
If $Gal(K/\Q) \cong S_4$, $$ d_M \le c q(\rho)^{3/2} $$ with an absolute constant $ c > 1$.
\end{myprop}

\begin{proof}
We use a strategy similar to the one we used earlier. Note that the representation $\rho$ of $S_4$ is the twist of the standard representation by the alternating character. This can be seen from the character table of $S_4$ since $\rho$ is irreducible and of trivial determinant.

{\bf Case 1: I is a cyclic subgroup of order 2.} There are two conjugacy classes, one containing 6 elements of order 2 and one containing 3 elements of order 2, which cover all the elements of order 2 in\S4. The character of our $3$-dimensional representation is valued $-1$ on both these classes. Hence, for restriction to any subgroup of order $2$, we see that $$ \mbox{Multiplicity of trivial character} = \frac{3(1) + (-1)(1)}{2} = 1. $$ Hence, $ \dim(V^I) = 1 $. So that $ \dim(V/V^I) = 2$.

{\bf Case 2: I is a cyclic subgroup of order 3.}  The character is valued $0$ on the unique conjugacy class of elements of order $3$. Hence, for restriction to any subgroup of order $3$, we see that $$ \mbox{Multiplicity of trivial character} = \frac{3(1) + 0(1) + 0(1)}{3} = 1. $$  Hence, $ \dim(V^I) = 1 $. So that $ \dim(V/V^I) = 2$.

{\bf Case 3: I is a cyclic subgroup of order 4.} Any subgroup of order $4$ contains, apart from the identity, two elements of order $4$ and one element of order $2$. Our character is valued $1$ on elements of order $4$ and $-1$ on elements of order $2$. (Elements of order $4$ all belong to the same conjugacy class and the class does not matter for elements of order $2$ as our character is valued the same on both of them as mentioned above.) Hence, for restriction to any subgroup of order $4$, we see that $$ \mbox{Multiplicity of trivial character} = \frac{3(1) + 1(1) + (-1)(1) + 1(1)}{4} = 1. $$ Hence, $ \dim(V^I) = 1 $. So that $ \dim(V/V^I) = 2$.

Therefore, by (\ref{eq:conductor}), we see that $$ q(\rho) \ge \prod_{\substack{p | q(\rho) \\ p > m}} p^2 $$ and combined with (\ref{eq:ram}), it completes the proof.

We remark that we can follow the alternate method mentioned in the previous proposition here as well.

\end{proof}

\begin{myprop} \label{SO3A5}
If $Gal(K/\Q) \cong A_5$, $$ d_M \le c q(\rho)^{2} $$ with an absolute constant $ c > 1$.
\end{myprop}

\begin{proof}
Since $m$ is now $5$ rather than $4$, we have a different bound \begin{equation}\label{eq:ram2} d_M \le c \prod_{\substack{p | d_K \\ p > 5}} p^4 \end{equation} with $ c = 2^{14}3^{9}5^{9}$ from the same reference \cite{Serre}.

Now, we again consider cases of cyclic subgroups. Note that, in this case, we have two $3$-dimensional representations.

{\bf Case 1: I is a cyclic subgroup of order 2.} All the elements of order $2$ in\Al5 are conjugate to each other. Both our characters are valued $-1$ on this class. Hence, for restriction of either of the representations to any cyclic subgroup of order $2$, we see that $$ \mbox{Multiplicity of trivial character} = \frac{3(1) + (-1)(1)}{2} = 1. $$  Hence, $ \dim(V^I) = 1 $. So that $ \dim(V/V^I) = 2$.

{\bf Case 2: I is a cyclic subgroup of order 3.} All the elements of order $3$ in\Al5 are conjugate to each other. Both our characters are valued $0$ on this class. Hence, for restriction of either of the representations to any cyclic subgroup of order $3$, we see that $$ \mbox{Multiplicity of trivial character} = \frac{3(1) + 0(1) + 0(1)}{3} = 1. $$  Hence, $ \dim(V^I) = 1 $. So that $ \dim(V/V^I) = 2$.

{\bf Case 2: I is a cyclic subgroup of order 5.} There are two conjugacy classes in\Al5, each containing 12 elements of order 5, which cover all the elements of order 5. As before, we can compute $ \dim(V^I)$ using a character table. But, in this case, it is more efficient to use the alternate method. We just note that since the determinant of $\rho$ is $1$, the image under $\rho$ of a non-trivial element of $I$ is conjugate to 

$$ 
\left( \begin{array}{ccc}
1 &   &   \\
  & \omega &  \\
  &   & \omega^2 \end{array} \right) $$ where $\omega$ is a primary $5$th root of unity. Thus, it is immediate that $ \dim(V^I) = 1 $. So that $ \dim(V/V^I) = 2$.

Therefore, by (\ref{eq:conductor}), we see that $$ q(\rho) \ge \prod_{\substack{p | q(\rho) \\ p > m}} p^2 $$ and combined with (\ref{eq:ram2}), it completes the proof.

\end{proof}

Finally, we have $$ \vartheta_1(x) = \sum_{ \substack{Gal(K / \Q) \subset  SO(3) \\ Gal(K / \Q) \cong A_4 \\ q(\rho) \le x }} 1 + \sum_{ \substack{Gal(K / \Q) \subset  SO(3) \\ Gal(K / \Q) \cong S_4 \\ q(\rho) \le x }} 1 + \sum_{ \substack{Gal(K / \Q) \subset  SO(3) \\ Gal(K / \Q) \cong A_5 \\ q(\rho) \le x }} 1,$$ which translates using propositions~\ref{SO3A4}, \ref{SO3S4}, \ref{SO3A5}, and recalling that $K$ is the normal closure of $M$, to $$ \vartheta_1(x) \le  \sum_{ \substack{ Gal(K / \Q) \cong A_4 \\ [M:\Q] = 4 \\ d_M \le cx^{\frac{3}{2}} }} 1 + \sum_{ \substack{ Gal(K / \Q) \cong S_4 \\ [M:\Q] = 4 \\ d_M \le cx^{\frac{3}{2}} }} 1 + \sum_{ \substack{ Gal(K / \Q) \cong A_5 \\ [M:\Q] = 5 \\ d_M \le cx^{2} }} 1. $$ Since \cite[Theorem $7$]{Bhargava01}, and \cite[Theorem 1]{BCT} which builds upon Theorem $6$ in \cite{Bhargava02}, show that $$\eta_{\Q, 4}^{A_4}(x) = O(x), \text{ } \eta_{\Q, 4}^{S_4}(x) = O(x), \text{ } \eta_{\Q, 5}^{A_5}(x) = O(x^{1- \beta}), $$ where $\beta$ is any positive constant less than $1/120$, we see that the above implies \begin{equation} \label{eq:firstterm} \vartheta_1(x) = O(x^{2 - 2 \beta}). \end{equation} This yields the required estimate for $\vartheta_{1}(x)$.

\end{section}

\begin{section}{Finite Subgroups of $ O(3) $ that are not \\ contained in $ SO(3) $ - Part $1$}
In this section and the next, we focus on bounding $\vartheta_2(x)$.

We deal with this case in two parts, depending upon whether $-1$, the negative of the identity matrix in $3$ dimensions, is in the image of $Gal(K / \Q)$. Part $1$ is devoted to the case :

 $$ Gal(K/ \Q) \not \subset SO(3), -1 \not \in Gal (K/ \Q). $$ 

In this case, we prove a lemma which straightaway tells us what $Gal(K/ \Q)$ is.

\begin{mylemma} Let $G$ be a finite irreducible subgroup of $O(3)$ which is not contained in $SO(3)$. Assume further that $-1 \not \in G$. Then $G \cong S_4$.

\end{mylemma}

\begin{proof}
Let $$ H = G \cap SO(3). $$ Then we can write $$ G = H \cup \kappa H, $$ where $ det(\kappa) = -1$.

 Then we can define$$  H^{*}  = H \cup (- \kappa)H, $$ so that $ H^{*} \subset SO(3) $. Note that $G$ and $H^{*}$ are ``isoclinic'', i.e. they only possibly differ by scalars. (Refer to \cite{CCNPW}, page xxiii, for more details on isoclinism.) Therefore, $H^{*}$ is an irreducible subgroup of \SO3. As a result, we see that $H^{*} \cong A_4, S_4 \text{ or } A_5$. But, \A4 or\Al5 do not possess index $2$ subgroups. Hence, $$ H^{*} \cong S_4. $$

Then, we can give an explicit isomorphism from $G$ to $H^{*}$ using the index $2$ subgroup $H$, viz. send $ h \rightarrow h $ and $ \kappa h \rightarrow (-\kappa)h $. It is easily seen that this is an isomorphism, and thus $ G \cong S_4 $.

\end{proof}

Indeed we see that\S4 has a faithful $3$-dimensional representation which satisfies our hypotheses. This is the standard representation mentioned in the introduction of this paper. Our earlier method of comparing conductors and discriminants yields weaker bounds than what are needed for our purpose, because we cannot rule out the possibility that $\det \rho$ is nontrivial on $I$. If $I$ is of order $2$, then the matrix corresponding to the non-trivial element of inertia would then be conjugate to $$ \left( \begin{array}{ccc}
-1 &   &   \\
  & 1 &  \\
  &   & 1 \end{array} \right) $$ which yields a weaker lower bound for $q(\rho)$ and thus allows for a larger number of $\rho$. This does not happen, but we need to take a different path to prove this, which we describe below.

The standard representation (i.e. the $3$-dimensional irreducible representation with non-trivial determinant) of\S4 is monomial. This fact can be inferred from the character table of $S_4$ and Frobenius reciprocity. In particular, it implies existence of a subfield $M$ of $K$ such that $[M: \Q] = 3$, $Gal(K/M) = D_8$, together with a quadratic $1$-dimensional character $\chi$ of $Gal(K/M)$ such that $ \rho = Ind_{M/ \Q} (\chi) $, where $Ind_{M / \Q}$ denotes induction of representation from the subgroup $Gal(K/M)$ to the group $Gal(K/ \Q)$. Now we count pairs $(M, \chi)$ to count representations $\rho$. We have classical results by Davenport and Heilbronn (\cite{DH}) on counting cubic number fields by discriminants. Since we are only concerned about the main term here, we do not need the error improvement by Bhargava, Shankar and Tsimerman (\cite{BST}). Note that for each such pair $(M, \chi)$, we have two more pairs $(M', \chi')$ and $(M'', \chi'')$, corresponding to conjugate copies of $D_8 \subset S_4$ that will give us the same representation, but that will only affect the constant term in our bounds.

We begin by denoting our counting function : $$ \Theta(x) := \sum_{\substack{\rho \\ Gal(K/ \Q) \not \subset SO_3 \\ -1 \not \in Gal(K/ \Q) \\ q(\rho) \le x}} 1. $$

The Conductor-Discriminant formula gives \begin{equation} \label{eq:cond} q(\rho) = d_Mq(\chi), \end{equation} where $d_M$ is the discriminant of the field $M$ and $q(\chi)$ is the absolute norm of the conductor of the character $\chi$. 

Thus, we see that  \begin{equation} \label{eq:ineq} \Theta(x) = \sum_{\substack{\rho \\ Gal(K/ \Q) \not \subset SO_3 \\ Gal(K/ \Q) \cong S_4 \\ q(\rho) \le x}} 1 \mbox{    } \le \sum_{\substack{(M, \chi) \\ [M: \Q] = 3 \\ \chi^2 = 1 \\ q(\chi)d_M \le x}} 1, \end{equation} where the first equality is due to the lemma.

We need to count the extensions $M$ as well as characters $\chi$. We write $\theta_{M,2}(x)$ for the number of characters $\chi$ of $M$ with $\chi^2 = 1$ and $q(\chi) \le x$. We shall use upper bounds for $\theta_{M,2}(x)$ from the appendix of this paper. These bounds are slightly weaker than the actual asymptotic if we work with a field $M$ that is fixed, but since we are working with varying fields $M$ at the same time, these bounds, which are uniform as long as the degree $[M:\Q]$ is fixed (which is true in our case), will work better, and the only expense incurred is a power of logarithm. It can be seen from the final proposition of this section below, that this increased power does not affect our result. We note that the asymptotic is an interesting result in itself, which follows from computing the residue of an appropriate Zeta Function and knowledge of bounds on the class number and the regulator of a number field.

From the corollary to Proposition $2$ in the appendix to this paper, we see that \begin{equation} \label{eq:bound} \theta_{M,2}(x) \ll \sqrt{d_M}(\mbox{log }d_M)^{2}x(\mbox{log }x)^{2}, \end{equation} where the implied constant is independent of $M$, since $c$ and $m$ are now fixed.

We now prove our main result of this section.

\begin{myprop} \label{O3S4}
Let $\rho, K$ be as before. Then \begin{equation*} \Theta(x) = O(x^{\frac{3}{2} + \epsilon}). \end{equation*}
\end{myprop}

\begin{proof}
By (\ref{eq:ineq}), it is sufficient to prove \begin{equation} \sum_{\substack{(M, \chi) \\ [M: \Q] = 3 \\ q(\chi)d_M \le x }} 1 = O(x^{\frac{3}{2} + \epsilon}), \end{equation} where $\chi$ is a quadratic character of $M$. That is, we wish to prove \begin{equation} \label{eq:charest} \sum_{\substack{M \\ [M: \Q] = 3 \\ d_M \le x}} \mbox{ } \sum_{\substack{\chi \\ q(\chi) \le x/d_M}} 1 =  O(x^{\frac{3}{2} + \epsilon}). \end{equation} From (\ref{eq:bound}), we see that \begin{equation} \sum_{\substack{M \\ [M: \Q] = 3 \\ d_M \le x}} \mbox{ } \sum_{\substack {\chi \\ q(\chi) \le x/d_M}} 1 \ll \sum_{\substack{[M: \Q] = 3 \\ d_M \le x }} \sqrt{d_M}(\mbox{log }d_M)^{2}(\mbox{log }\frac{x}{d_M})^{2} \frac{x}{d_M}. \end{equation} The implied constant is uniform. Hence, we get \begin{equation} \sum_{\substack{M \\ [M: \Q] = 3 \\ d_M \le x}} \mbox{ } \sum_{\substack{\chi \\ q(\chi) \le x/d_M}} 1 = x(\mbox{log }x)^{2} O \left( \sum_{\substack{M \\ [M: \Q] = 3 \\ d_M \le x} } \frac{(\mbox{log } d_M)^{2}}{\sqrt{d_M}} \right). \end{equation} \cite[Theorem $1$]{DH} yields \begin{equation} \sum_{\substack{M \\ [M: \Q] = 3 \\ d_M \le x} } 1 \sim cx, \end{equation} where $c$ is an absolute constant. Hence, the above sum can be estimated as \begin{equation} \left( \sum_{\substack{M \\ [M: \Q] = 3 \\ d_M \le x} } \frac{(\mbox{log } d_M)^{2}}{\sqrt{d_M}} \right) = O(x^{\frac{1}{2} + \epsilon}), \end{equation} which proves (\ref{eq:charest}). Thus, the proposition follows.

\end{proof}

\end{section}

\begin{section}{Finite Subgroups of $ O(3) $ that are not \\ contained in $ SO(3) $ - Part $2$}
We deal with the remaining cases in this section. These cases are characterized by : $$ Gal(K/ \Q) \not \subset SO(3), -1 \in Gal (K/ \Q). $$ Put $H = Gal(K/ \Q) \cap SO(3) $. Then we see that $ Gal(K/ \Q) \cong H \times \{ \pm 1 \} $. The Artin representation $\rho$ we are considering can be written as $ \rho \cong \sigma \bigotimes \epsilon $, where $ \sigma $ is the irreducible three-dimensional representation of $H$ given by the inclusion $H \subset SO(3)$ and $ \epsilon $ is a quadratic character of $Gal(\overline{\Q}/ \Q)$. Thus, we now estimate the number of pairs $ (\sigma , \epsilon) $ and obtain the bounds that we need.

We first do this in the case where $ H \cong $\Al5. Let us denote the corresponding\Al5-subextension of $K$ by $L$. In this case, the representation $ \sigma $ of\Al5 is a primitive representation, since $A_5$ has no index three subgroups, and $\sigma$, if at all induced from a proper subgroup, would have to be induced from such a subgroup in order to be $3$-dimensional. (That $A_5$ does not possess such a subgroup can be seen from the fact that the existence of such a subgroup would yield a nontrivial homomorphism from $A_5$ to $S_3$ due to the action on coset space, kernel of which would be a nontrivial proper normal subgroup of $A_5$.) In fact, there are two representations of\Al5 possible, and both of them are primitive. As before, let $M$ be a subfield of $L$ such that $[M: \Q]$ = 5. The normal closure of $M$ is $L$ for any choice of $M$.

We wish to estimate $$ \Psi(x) : = \sum_{\substack{\rho \\ q(\rho) \le x \\ \rho \cong \sigma \bigotimes \epsilon}} 1, $$ where $\sigma$ is a faithful irreducible representation of $A_5 \cong Gal(L/ \Q)$, considered as an Artin representation of $\Q$, and $\epsilon$ is a quadratic character of $Gal(\overline{\Q} / \Q)$.

For a fixed $ \sigma $ we look at $ q(\sigma \otimes \epsilon) \mbox{ and } q(\sigma)$.  Let $$ q_{tame}(\sigma) = \prod_{p \in X} p^{e_p} $$ be the tame conductor of $\sigma$, where $X$ is the set of tamely ramified primes in $L$. By a computation done previously, each $e_p = 2$. Thus, \begin{equation} \label{eq:tame} q_{tame}(\sigma) = \prod_{p \in X} p^{2}. \end{equation} Let $X = A \cup B \cup C$, where $A, B,$ and $C$ are the sets of tamely ramified primes with inertia subgroup isomorphic to $\Z / 2 \Z, \Z / 3 \Z,$ and $\Z / 5 \Z$ respectively. 

Let us write the conductor of $\epsilon$ as \begin{equation}\label{eq:condepsilon} q(\epsilon) = 2^{\alpha} \prod_{j = 1}^{s} q_{j}, \end{equation} where $\alpha = 0, 2,$ or $3,$ and $q_j$ are distinct odd primes. We can write $ \epsilon = \chi \chi' $ where $\chi$ and $\chi'$ are quadratic characters of $\Q$ such that $$ p | q(\chi) \iff p =2,3,5 \mbox{ or } p \in X. $$ Thus, primes that divide $q(\chi')$ are  primes that are unramified in $L$ which divide the conductor of $\epsilon$.

\begin{myprop}
Let $\sigma, \epsilon, \chi, \chi'$ be as above. Then $$ q(\sigma \otimes \epsilon) = q(\sigma \otimes \chi \chi') = q(\sigma \otimes \chi) q(\chi')^{3}. $$
\end{myprop}

\begin{proof}
This follows from a well-known local computation for each prime dividing $q(\sigma \otimes \chi)$ and $q(\chi')$. See \cite[Corollary $1^{\prime}$]{SerLoc} for details.  
\end{proof}

Thus, we see that the condition $$ q(\rho) = q(\sigma \otimes \epsilon) \le x $$ translates to the condition \begin{equation} \label{eq:translate} q(\sigma \otimes \chi) \le x \mbox{ and } q(\chi') \le \left( \frac{x}{q(\sigma \otimes \chi)} \right)^{\frac{1}{3}}. \end{equation}

Hence, we get \begin{equation} \label{eq:unrub} \Psi(x) \le \sum_{q(\sigma \otimes \chi) \le x} \mbox{ } \sum_{q(\chi') \le \left( \frac{x}{q(\sigma \otimes \chi)} \right)^{\frac{1}{3}} }  1. \end{equation}

The number of quadratic characters of $\Q$ with conductor $\le x$ is $O(x)$. So (\ref{eq:unrub}) yields $$ \Psi(x) \ll \sum_{q(\sigma \otimes \chi) \le x} \left( \frac{x^{\frac{1}{3}}}{ q(\sigma \otimes \chi)^{\frac{1}{3}}} \right), $$ which gives the following proposition.

\begin{myprop}

\begin{equation} \label{eq:ub} \Psi(x) \ll x^{\frac{1}{3}} \sum_{q(\sigma \otimes \chi) \le x} \frac{1}{q(\sigma \otimes \chi)^{\frac{1}{3}}}. \end{equation}

\end{myprop}

Now, let $$\Theta(x) = \sum_{q(\sigma \otimes \chi) \le x} 1, $$ where $\sigma$ and $\chi$ are as above.

\begin{myprop}

We have $$ \Theta(x) = O(x^{2 - 2\beta}), $$ where $\beta$ is any positive constant less than $\frac{1}{120}$.

\end{myprop}

\begin{proof}
We convert the problem of estimating $ \Theta(x)$ into a problem of counting $A_5$-extensions of $\Q$ and quadratic characters of $\Q$. For this, we look at the conductors and discriminants.

Since we only let $2,3,5$ or primes that are tamely ramified in $L$ remain in the conductor of $\chi$, we see that $$ q(\chi) = 2^{\alpha} 3^{\beta} 5^{\gamma} \prod_{p \in Y} p $$ for some $Y \subset X$ and $ \beta, \gamma \in \{ 0, 1 \} $. 

Let $A' = A \cap Y, B' = B \cap Y, C' = C \cap Y$ and $A'' = A \setminus A', B'' = B \setminus B', C'' = C \setminus C'$.

We can compute the effect of twisting by $\chi$ at each prime locally by looking at image under $\rho$ of $I$, {\it cf.} proof of Proposition~\ref{SO3A4}. We see using equations (\ref{eq:tame}) and (\ref{eq:condepsilon}) that the tame conductor of $\sigma \otimes \chi$ is given by \begin{equation} \label{eq:twistcond} q_{tame}(\sigma \otimes \chi) = \left( \prod_{p \in A'} p \prod_{p \in A''} p^{2} \right) \left( \prod_{p \in B'} p^{3} \prod_{p \in B''} p^{2} \right) \left( \prod_{p \in C'} p^{3} \prod_{p \in C''} p^{2} \right). \end{equation}

On the other hand, by a computation involving ramification degrees, we obtain a bound on the tame discriminant of $M$ : \begin{equation} \label{eq:tamedisc} d^{tame}_{M} \le \prod_{p \in A} p^{2} \prod_{p \in B} p^{2} \prod_{p \in C} p^{4}. \end{equation}

Comparing the above expressions, we see that \begin{equation} \label{eq:tamediscbound} d^{tame}_{M} \le q_{tame}(\sigma \otimes \chi)^{2}. \end{equation}

This is the inequality which helps us translate the bound on the conductor to a bound on the discriminants. We have not yet dealt with the primes $2,3,5$ which might be wild primes. But we have a uniform bound for them as we have seen before, cf. \cite{Serre}. Letting $c = 2^{14}3^{9}5^{9}$, we see from \emph{loc. cit.} and (\ref{eq:translate}) that \begin{equation} \label{eq:discbound} d_M \le c d^{tame}_{M} \le c q_{tame}(\sigma \otimes \chi)^{2} \le c q(\sigma \otimes \chi)^{2} \le cx^{2}. \end{equation}

Equipped with these results, we now obtain our result.  We have $$ \Theta(x) \le \sum_{\sigma} \sum_{\chi} 1 $$ where $\sigma$ and $\chi$ are as above and the sum runs only over pairs such that $q(\sigma \otimes \chi) \le x$. Thus, we get \begin{equation} \Theta(x) \ll \sum_{\substack{M \\ d_M \le cx^2}} \sum_{\substack{\chi \\ \chi^{2} = 1 \\ q(\chi) | 8d_M}} 1. \end{equation} 

The inner sum is over quadratic characters $\chi$ and is thus $O(x^{\epsilon})$ once $\sigma$ is fixed, because the number of divisors of $d_M$ is $O(d_{M}^{\epsilon})$. Hence, we get \begin{equation} \Theta(x) \ll x^{\epsilon} \sum_{\substack{M \\ d_M \le cx^2}} 1. \end{equation} We then appeal to \cite{BCT} as before, to obtain \begin{equation} \Theta(x) \ll x^{\epsilon} \left( x^{2} \right)^{1- \beta} \end{equation} where $\beta$ is a positive constant less than $1/120$, which finishes our proof.

\end{proof}

Using this proposition, we obtain our main result :

\begin{myprop} \label{O3A5}
\begin{equation} \label{eq:primO3} \Psi(x) = O(x^{2 - 2\beta}) \end{equation} where $\beta$ is as before.
\end{myprop}

\begin{proof}
This follows from combining the previous propositions $6$ and $7$ using Abel partial summation.
\end{proof}

We move on to the remaining cases. Recall that $ Gal(K / \Q) \cong H \times \{ \pm 1 \}$, where $H$ is an irreducible finite subgroup of $SO(3)$. We have dealt with the case $ H \cong $\Al5. We deal with the cases $ H \cong $ \A4, $ H \cong$\S4 below. Let $L/ \Q$ be a subextension such that $Gal(L/ \Q) \cong H$.

In these cases, the representations $\rho$ can be again written as $ \sigma \otimes \epsilon $, where $\sigma$ is an irreducible $3$-dimensional representation of \A4 or\S4 with trivial determinant. Such representations $\sigma$ are necessarily monomial, say induced from a cubic subextension $M/ \Q$ of $L$. As before, this follows from looking at the character tables of respective groups. $\epsilon$ is a quadratic character of $Gal(\overline{\Q}/ \Q)$ as before. Note that $\epsilon = \det \rho$.

\begin{myprop} \label{O3A4S4}
Let $K, M, H$ be as before where $H \cong $ \A4 or\S4. Define $$ \Phi(x) := \sum_{\substack{\rho \\ q(\rho) \le x \\ \rho \cong \sigma \otimes \epsilon}} 1 $$ where $\rho, \sigma, \epsilon$ are also as before.

Then, \begin{equation} \label{eq:monoms4} \Phi(x) = O(x^{\frac{3}{2} + \epsilon'}) \end{equation} where $\epsilon'$ is arbitrarily small.
\end{myprop}

\begin{proof}
We prove the result for $H \cong $\S4. The case $ H \cong  A_4$ is completely analogous. Let $Ind_{M/ \Q}$ and $res_{M/ \Q}$ denote the induction and restriction functors for the group $Gal(L/\Q)$ and subgroup $Gal(L/M)$. If $\sigma = Ind_{M/ \Q}(\chi)$ for some $\chi$, we see, by Frobenius reciprocity, that $$ \rho \cong Ind_{M/ \Q}(\chi \otimes res_{M/ \Q} \epsilon). $$ Here, $M$ determines $\sigma$ which in turn fixes $\epsilon$. Moreover, both $\chi$ and $\epsilon$ are quadratic characters. We denote the quadratic character $\chi \otimes \text{res } \epsilon$ by $\chi'$.

This is analogous to our methods in section $5$, where we counted pairs $(M, \chi')$ where $M$ is a cubic extension of $\Q$ and $\chi'$ is a quadratic character of $M$. We again write : $$ q(\rho) = d_M q(\chi'). $$ Then, using the same method, from corollary to Proposition 2 in the appendix, we have : $$ \theta_{M, 2}(x) \ll \sqrt{d_M}(\mbox{log }d_M)^{2}x(\mbox{log }x)^{2}. $$ And thus, we get $$ \Phi(x) \ll \sum_{\substack{[M: \Q] = 3 \\ d_M \le x }} \sqrt{d_M}(\mbox{log }d_M)^{2}(\mbox{log }\frac{x}{d_M})^{2} \frac{x}{d_M}, $$ from which we get the required result.

\end{proof}

\end{section}

\begin{section}{Proof of Theorem 1}
By Propositions~\ref{O3S4}, \ref{O3A5}, \ref{O3A4S4}, we see that \begin{equation} \label{eq:secondterm} \vartheta_2(x) = O(x^{2 - 2 \beta}) \end{equation} where $\beta$ is any positive constant less than $1/120$, which coupled with equations (\ref{eq:division}) and (\ref{eq:firstterm}) finishes the proof.
\end{section}


\begin{thebibliography}{9}

\bibitem{ABCZ} G. Anderson, D. Blasius, R. Coleman, and G. Zettler,
\emph{On representations of the Weil group with bounded conductor,}
Forum Math. 6 (1994), 537 -- 545.

\bibitem{Artin}M. Artin, 
\emph{Algebra, Second Edition,}
Pearson, 2 edition (August 13, 2010).

\bibitem{Bhargava01} M. Bhargava,
\emph{The density of discriminants of quartic rings and fields,}
Ann. Math. 162 (2005), 1031 -- 1063.

\bibitem{Bhargava02} M. Bhargava,
\emph{The density of discriminants of quintic rings and fields,}
Ann. Math. 172 (2010), 1559 -- 1591.

\bibitem{BCT} M. Bhargava, A. Cojocaru, and F. Thorne,
\emph{The square sieve and the number of $A_5$-quintic extensions of bounded discriminant,}
To appear.

\bibitem{BG} M. Bhargava, E. Ghate,
\emph{On the average number of octahedral newforms of prime level,}
Math. Ann. 344 (2009), no. 4, 749768.

\bibitem{BST} M. Bhargava, A. Shankar, J. Tsimerman,
\emph{On the Davenport-Heilbronn theorems and second order terms,}
Invent. Math. 193(2) (2013), 439 -- 499.

\bibitem{CCNPW} J. Conway, R. Curtis, S. Norton, R. Parker, R. Wilson,
\emph{Atlas of finite groups: maximal subgroups and ordinary characters for simple groups,}
Oxford, England: Clarendon Press, (1985).

\bibitem{DH} H. Davenport and H. Heilbronn,
\emph{On the density of discriminants of cubic fields II,}
Proc. Roy. Soc. London Ser. A 322 (1971), no. 1551, 405 -- 420.

\bibitem{FH} W. Fulton, J. Harris,
\emph{Representation Theory : A First Course (Readings in Mathematics),}
Springer, Corrected edition (October 22, 1991).

\bibitem{K} J. Kl\"{u}ners,
\emph{The number of $S_4$-fields with given discriminant,}
Acta Arith. 122 (2006), no. 2, 185 -- 194.

\bibitem{Rohrlich} D. Rohrlich,
\emph{Self-dual Artin representations,}
Automorphic Representation and L-functions, D. Prasad, C.S. Rajan, A. Sankaranarayanan, J. Sengupta, eds., Tata Institute of Fundamental Research Studies in Mathematics Vol. 22 (2013), 455 -- 499. 

\bibitem{SerLoc} J.-P. Serre,
\emph{Local Fields : Graduate Texts in Mathematics {\bf 67},}
Springer, Second corrected edition (1995).

\bibitem{Serre} J.-P. Serre,
\emph{Quelques applications du th\'{e}or\`{e}me de densit\'{e} de Chebotarev,}
Publ. Math. IHES 54 (1981), 123 -- 201. (=\emph{Oeuvres} vol. III, no. 125.)

\bibitem{T} M. J. Taylor,
\emph{On the equidistribution of Frobenius in cyclic extensions of a number field,}
J. London Math. Soc. 29 (1984), 211 -- 213.


\end{thebibliography}
\end{document}

\centerline
{{\bf Appendix:}}

\centerline
{{\bf Ray class characters of bounded order and bounded conductor}}
\vskip .2 cm

\centerline {{David E. Rohrlich}}
\vskip .1 cm

We introduce a partial order on 
the set of formal
Dirichlet series
with nonnegative real coefficients.
Given two such series $A(s)=\sum_{q\geqslant 1} a(q) q^{-s}$
and
$B(s)=\sum_{q\geqslant 1} b(q) q^{-s}$, 
write
$A(s) \preccurlyeq B(s)$ to mean that
$a(q) \leqslant b(q)$ for all $q\geqslant 1$.
It is readily verified that if 
$A(s) \preccurlyeq B(s)$ and $C(s) \preccurlyeq D(s)$
then $A(s)C(s) \preccurlyeq B(s)D(s)$.
Furthermore the implication holds at the level of
Euler products: if $A(s)=\prod_p A_p(s)$
and $B(s)=\prod_pB_p(s)$ with $A_p(s)\preccurlyeq B_p(s)$
for all $p$ then  
$A(s)\preccurlyeq B(s)$.

By way of illustration, let $M$ be a number field, $\mathcal O_M$
its ring of integers, and
$\zeta_M(s)$ the associated
Dedekind zeta function. Then it is a standard remark that 
\begin{equation}\label{zetakzeta}
\zeta_M(s) \preccurlyeq \zeta(s)^{m},
\end{equation}
where $m=[M:\Q]$.
Indeed let $p$ be a rational prime
and $\mathfrak p$ a prime ideal of $\mathcal O_M$ above
$p$, say of residue class degree $f$. The Euler factor of $\zeta_M(s)$
at $\mathfrak p$ satisfies
\[
(1-(\mathbf N\mathfrak p)^{-s})^{-1}=\sum_{\nu\geqslant 0}p^{-\nu fs}
\preccurlyeq \sum_{\nu\geqslant 0}p^{-\nu s}=(1-p^{-s})^{-1}.
\]
Hence if there are exactly $r$ prime ideals $\mathfrak p$ above $p$ then
\[
\prod_{\mathfrak p |p}(1-(\mathbf N\mathfrak p)^{-s})^{-1}
\preccurlyeq (1-p^{-s})^{-r}
\preccurlyeq (1-p^{-s})^{-m}.
\]
Passing to Euler products we obtain \eqref{zetakzeta}.

It follows from the definitions that if $A(s)\preccurlyeq B(s)$
then the associated summatory functions
$\vartheta_A(x)=\sum_{n\leqslant x} a(n)$
and 
$\vartheta_B(x)=\sum_{n\leqslant x} b(n)$
satisfy $\vartheta_A(x)\leqslant \vartheta_B(x)$
for all $x$.
For example, let $A(s)$ and
$B(s)$ be the two sides of \eqref{zetakzeta}: 
Using Theorem 7.7 on p.\ 154 of \cite{bateman} to 
estimate 
the summatory function
of $\zeta(s)^{m}$, we obtain 
\begin{equation}\label{hereitis}
\sum_{\mathbf N\mathfrak q \leqslant x} 1 \ll x(\log x)^{m-1}
\phantom{................} (x\geqslant 2),
\end{equation}
where $\mathfrak q$ denotes a nonzero ideal of $\mathcal O_M$ and
the implicit constant depends only on $m$, not on $M$.

To illustrate the use of \eqref{hereitis}, let us deduce a standard bound
for the class number $h_M$ of $M$. Let
$r_1$ and $r_2$ be 
the number of real embeddings and half the number
of complex
embeddings of $M$, so that $r_1+2r_2=m$.
Thus the Minkowski constant $(4/\pi)^{r_2} m!/m^m$ is
bounded above by
\[
\mink=(4/\pi)^{m/2}\dfrac{m!}{m^m},
\]
and therefore Minkowski's theorem gives
\begin{equation}\label{minkowski}
h_M \leqslant \sum_{\mathbf N\mathfrak q \leqslant \mink\sqrt{d_M}} 1,
\end{equation}
where $d_M$ is the absolute value of the discriminant of $M$
(cf.\ \cite{lang}, pp.\ 119-120).
Combining \eqref{minkowski} with
\eqref{hereitis}, we recover the well-known bound
\begin{equation}\label{hk}
h_M\ll \sqrt{d_M}(\log d_M)^{m-1}
\phantom{................} (\mink\sqrt{d_M} \geqslant 2),
\end{equation}
where the implicit constant depends only on $m$.
We shall regard $m$ as a fixed integer $\geqslant 2$, and thus the condition
$\mink\sqrt{d_M} \geqslant 2$ is satisfied for all but finitely many $d_M$
with $[M:\Q]=m$. Furthermore, since $m\geqslant 2$, we have $d_M \geqslant 2$.
Therefore we can remove the condition 
$\mink\sqrt{d_M} \geqslant 2$ from \eqref{hk} and still 
assert that the implicit constant in \eqref{hk} depends only on $m$.
Actually it is more useful to state \eqref{hk} for
$h_M^\nar$, the narrow ray class number 
of $M$. Since $h_M^{\nar}\leqslant 2^{r_1}h_M$,
we have
\begin{equation}\label{hkbis}
h_M^{\nar}\ll \sqrt{d_M}(\log d_M)^{m-1},
\end{equation}
where the implicit constant depends only on $m$.

It is convenient to refine
the relation $\preccurlyeq$ slightly.
Suppose that $A(s)$ and $B(s)$
are Dirichlet series over $M$ in the sense that they 
are presented to us in the form $A(s)=\sum_{\mathfrak q} a(\mathfrak q)
(\mathbf N \mathfrak q)^{-s}$
and $B(s)=\sum_{\mathfrak q} b(\mathfrak q)
(\mathbf N \mathfrak q)^{-s}$, where $\mathfrak q$ denotes as before
a nonzero ideal of $\mathcal O_M$.
We write
$A(s) \preccurlyeq _M B(s)$
to mean 
that 
$a(\mathfrak q)\leqslant b(\mathfrak q)$
for all $\mathfrak q$. 
Thus $\preccurlyeq$ coincides with $\preccurlyeq _\Q$.
Of course every Dirichlet series is a Dirichlet series over
$\Q$, and one readily verifies that if
$A(s) \preccurlyeq_M
 B(s)$ then
$A(s) \preccurlyeq B(s)$.
%Indeed when $A(s)$ and $B(s)$ are written as Dirichlet series
%over $\Q$, the coefficient of $n^{-s}$ is
%$\alpha(n)=\sum_{\mathbf N\mathfrak n =n} a(\mathfrak n)$
%and $\beta(n)=\sum_{\mathbf N\mathfrak n =n} b(\mathfrak n)$
%respectively. Hence $\alpha(n)\leqslant \beta(n)$
%for all $n$ if $a(\mathfrak n)\leqslant b(\mathfrak n)$
%for all $\mathfrak n$.

Given
a rational integer $c\geqslant 2$, let
\[ 
R_{M,c}(s)=
\sum_{\mathfrak q} h_{M,c}^*(\mathfrak q) (\mathbf N\mathfrak q)^{-s}
\]
where
$h_{M,c}^*(\mathfrak q)$
is the number of idele class characters $\chi$ of
$M$ of conductor $\mathfrak q$
such that $\chi^c=1$.
Also put
\[
E_{M,c}(s)=
\prod_{p|c}\prod_{\mathfrak p|p}
\big(
\sum_{\nu=0}^{e(\mathfrak p)(v_p(c)+1)}(\mathbf N\mathfrak p)^{\nu(1-s)}\big),
\]
where $e(\mathfrak p)$ is the 
ramification index of $\mathfrak p$ 
over $p$ and
$v_p(c)$ the $p$-adic valuation of $c$.

\begin{proposition}\label{drayclass} 
$
R_{M,c}(s) \preccurlyeq _M
h_M^{\nar}\cdot (\zeta_M(s)/\zeta_M(2s))^{c-1}
\cdot E_{M,c}(s).
$
\end{proposition}

Define 
\[
E_{m,c}= \prod_{p|c}\prod_{e=1}^m\prod_{f=1}^m
\big(\sum_{\nu=0}^{e(v_p(c)+1)} p^{f\nu(1-s)}\big)^m,
\]
The following variant of Proposition 
\ref{drayclass} is weaker but actually more
useful:

\begin{proposition} \label{variant}
$
R_{M,c}(s) \preccurlyeq 
h_M^{\nar}\cdot \zeta(s)^{m(c-1)}
\cdot E_{m,c}(s).
$
\end{proposition}

\begin{proof} By inspection, 
$
E_{M,c}(s)\preccurlyeq E_{m,c}(s).
$
Also
\[
\zeta_M(s)/\zeta_M(2s) = \prod_{\mathfrak p}
(1+(\mathbf N\mathfrak p)^{-s})
\preccurlyeq \prod_p(\sum_{\nu\geqslant 0}p^{-\nu s})^m =\zeta(s)^m,
\]
where $\mathfrak p$ runs over all nonzero prime ideals of $\mathcal O_M$.
\end{proof}

Let $\vartheta_{M,c}(x)$ and $\vartheta_{m,c}(x)$
denote the summatory function
associated to $R_{M,c}(s)$ and 
$\zeta(s)^{m(c-1)}
\cdot E_{m,c}(s)$ respectively.
Then Proposition 
\ref{variant} gives 
\[
\vartheta_{M,c}(x)\leqslant
h_M^{\nar}\vartheta_{m,c}(x),
\]
which in conjunction with 
\eqref{hkbis} becomes
\begin{equation}\label{pinnacle}
\vartheta_{M,c}(x)\ll \sqrt{d_M}(\log d_M)^{m-1}\vartheta_{m,c}(x).
\end{equation}
Here the implicit constant depends only on $m$.
Since $E_{m,c}(s)$ is entire while $\zeta(s)$
has a simple pole at $s=1$, 
we obtain (cf.\ \cite{bateman}, {\it loc. cit.}):

\begin{corollary}
$\vartheta_{M,c}(x)\ll 
\sqrt{d_M}(\log d_M)^{m-1}
x(\log x)^{m(c-1)-1}$,
the implicit constant depending only on $c$ and $m=[M:\Q]$.
\end{corollary}

We turn to the proof of Proposition \ref{drayclass}. Put
\[
\varphi_M(\mathfrak q)=|(\mathcal O_M/\mathfrak q)^{\times}|,
\]
and let $\A^{\times}_M$ be the group of ideles of $M$. 
As usual, we think of $\A^{\times}_M$ as the 
restricted direct product ${\prod}'_v M_v^{\times}$,
where $v$ runs
over the places of $M$ and $M_v$ is the completion of
$M$ at $v$, and we identify
$M^{\times}$ 
with its image in $\A_M^{\times}$
under the diagonal embedding. We also put
\begin{equation}\label{adelicunits}
{\widehat{\mathcal O}}_M=\prod_{v\nmid \infty}\mathcal O_v,
\end{equation}
where $v$ runs over the finite places of $M$ and
$\mathcal O_v$ is the ring of integers of $M_v$.
By appending the coordinate 1 at the infinite places,
we may view
${\widehat{\mathcal O}_M}^{\times}$ as a subgroup of $\A_M^{\times}$.
Similarly, the product
$M_{\infty}^{\times}=\prod_{v|\infty} M_v^{\times}$
and its identity component
$(M_{\infty}^{\times})^0$ are subgroups of 
$\A_M^{\times}$ with coordinate 1 at the finite places.
With these conventions,
\[
h_M^{\nar} =
|\A_M^{\times}/(M^{\times}\cdot {\widehat{\mathcal O}_M}^{\times}
\cdot (M_{\infty}^{\times})^0)|
\]
(cf.\ \cite{lang}, pp.\ 146-147).
As idele class characters are trivial on 
the principal ideles and idele class characters
of finite order are 
trivial on the identity component at infinity,
we deduce that there are at most $h_M^{\nar}$ extensions
of a given character of ${\widehat{\mathcal O}_M}^{\times}$ 
to a finite-order idele class character
of $M$. 
Let us write 
$\varphi_{M,c}^{*}(\mathfrak q)$ for the number of characters $\chi$
of ${\widehat{\mathcal O}_M}^{\times}$ 
of order dividing $c$ and conductor $\mathfrak q$,
the conductor of a character of ${\widehat{\mathcal O}_M}^{\times}$
being defined in the same way as for idele class
characters. Then the preceding discussion gives
\[
h^*_{M,c}(\mathfrak q)\leqslant h_M^{\nar} 
\varphi_{M,c}^{*}(\mathfrak q).
\]
Now $\varphi_{M,c}^{*}$ is multiplicative because
${\widehat{\mathcal O}}^{\times}_M=\prod_{v\nmid \infty}\mathcal O^{\times}_v$
by
\eqref{adelicunits}.
Thus
\begin{equation}\label{atlast}
\sum_{\mathfrak q}
h^*_{M,c}(\mathfrak q) (\mathbf N\mathfrak q)^{-s} \preccurlyeq _M
h_M^{\nar}
\prod_{\mathfrak p}(\sum_{\nu\geqslant 0}\varphi_{M,c}^{*}(\mathfrak
p^{\nu})(\mathbf N\mathfrak p)^{-\nu s}),
\end{equation}
where $\mathfrak p$ runs over the nonzero prime ideals of $\mathcal O_M$.

We now focus on the Euler factor in \eqref{atlast}
corresponding to a particular prime ideal $\mathfrak p$.
Let $v$ be the corresponding place of $M$
and
$p$ the residue characteristic of $\mathfrak p$.
We consider cases according as $p|c$
or $p\nmid c$. 
In both cases we use the fact 
that if $\nu \geqslant 2$ then 
$\varphi_{K,c}^{*}(\mathfrak p^{\nu})$
is the number of characters of
$\mathcal O_v^{\times}$ 
of order dividing $c$ which
factor through
$\mathcal O_v^{\times}/(1+\mathfrak p^{\nu}\mathcal O_v)$ 
but not through
$\mathcal O_v^{\times}/(1+\mathfrak p^{\nu-1}\mathcal O_v)$.

Suppose first that
$p\nmid c$. Then any character of $\mathcal O_v^{\times}$
of order dividing $c$ is trivial on 
the pro-$p$-group
$1+\mathfrak p\mathcal O_v$. Hence if $\nu\geqslant 2$ then
$\varphi_{M,c}^{*}(\mathfrak p^{\nu})=0$. Furthermore
\[
\varphi_{M,c}^{*}(\mathfrak p)=\gcd (c, \mathbf N\mathfrak p -1)-1
\]
because $\mathcal O_v^{\times}/(1+\mathfrak p\mathcal O_v)$ is cyclic
and the trivial character of
$\mathcal O_v^{\times}$
does not have conductor
$\mathfrak p$.
In particular we have 
$
\varphi_{M,c}^{*}(\mathfrak p)\leqslant c-1,
$
whence
\[
\sum_{\nu\geqslant 0}\varphi_{M,c}^{*}(\mathfrak
p^{\nu})(\mathbf N\mathfrak p)^{-\nu s}
\preccurlyeq_M 1+ (c-1) (\mathbf N\mathfrak p)^{-s}.
\]
Therefore
\begin{equation}\label{binomial}
\sum_{\nu\geqslant 0}\varphi_{M,c}^{*}(\mathfrak
p^{\nu})(\mathbf N\mathfrak p)^{-\nu s}
\preccurlyeq _M (1+ (\mathbf N\mathfrak p)^{-s})^{c-1}
\end{equation}
by the binomial theorem.

Next suppose that $p|c$. If
$k \geqslant e(\mathfrak p)/(p-1)+1$ 
then every element
of $1+c\mathfrak p^{k} \mathcal O_v$
is a $c$th power (cf.\ \cite{lang}, p.\ 186). 
In particular,
every element of $1+c\mathfrak p^{e(\mathfrak p) +1} \mathcal O_v$
is a $c$th power. It follows that
$\varphi_{M,c}^{*}(\mathfrak p^{\nu}\mathcal O_v)=0$
for $\nu \geqslant e(\mathfrak p) (v_p(c)+1)+1$. 
Now 
for $1\leqslant \nu \leqslant 
e(\mathfrak p)(v_p(c)+1)$ we apply the trivial estimate
\[
\varphi_{M,c}^{*}(\mathfrak
p^{\nu})\leqslant |\mathcal O_v^{\times}/(1+\mathfrak p^{\nu}
\mathcal O_v)|.
\]
Since
$|\mathcal O_v^{\times}/(1+\mathfrak p^{\nu}\mathcal O_v)|=
(\mathbf N \mathfrak p)^{\nu-1}
(\mathbf N \mathfrak p -1)\leqslant
(\mathbf N \mathfrak p)^{\nu},
$
we obtain
\begin{equation}\label{nonbinomial}
\sum_{\nu\geqslant 0}\varphi_{M,c}^{*}(\mathfrak
p^{\nu})(\mathbf N\mathfrak p)^{-\nu s}
\preccurlyeq_M 
\sum_{\nu=0}^{ e(\mathfrak p)(v_p(c)+1)}(\mathbf N\mathfrak p)^{\nu(1-s)}.
\end{equation}
This completes our discussion of the individual Euler factors in 
\eqref{atlast}.

Now combine \eqref{atlast},
\eqref{binomial}, and \eqref{nonbinomial}. We obtain
\begin{equation}\label{almostthere}
\sum_{\mathfrak q}
h_{M,c}^*(\mathfrak q) (\mathbf N\mathfrak q)^{-s} \preccurlyeq_M
h_M^{\nar}\cdot \prod_{p\nmid c}\prod_{\mathfrak p|p}
(1+ (\mathbf N\mathfrak p)^{-s})^{c-1}
\cdot
E_{M,c}(s)
\end{equation}
We may weaken the estimate 
in \eqref{almostthere}
by extending the product over $p\nmid c$ 
to a product over all $p$, and then we use the identity
\[ 
\zeta_M(s)/\zeta_M(2s)=
\prod_{\mathfrak p} (1+ (\mathbf N\mathfrak p)^{-s}).
\]
Making this substitution  
in \eqref{almostthere},
we obtain Proposition \ref{drayclass}.